\documentclass[a4paper]{article}

\usepackage{amsmath}
\usepackage{amssymb}
\usepackage{amsthm}
\usepackage{indentfirst}

\makeatletter

\newcommand{\Rmnum}[1]{\uppercase\expandafter{\romannumeral #1}}

\begin{document}

\title{\textbf{Proof of Murphy-Cohen Conjecture on One-dimensional
Hard Ball Systems}}
\author{Lizhou Chen\footnotemark}
\date{}
\maketitle

\renewcommand{\thefootnote}{*}
\footnotetext{Institute of Mathematics, Fudan University, Shanghai
200433, China.\\\hspace*{15pt}E-mail: 031018007@fudan.edu.cn}

\noindent{\textbf{Abstract}} {\small{\ We prove the Murphy and
Cohen's conjecture that the maximum number of collisions of $n+1$
elastic particles moving freely on a line is $n(n+1)\over2$ if no
interior particle has mass less than the arithmetic mean of the
masses of its immediate neighbors. In fact, we prove the stronger
result that, for the same conclusion, the condition no interior
particle has mass less than the geometric mean, rather than the
arithmetic mean, of the masses of its immediate neighbors
suffices.}



\section{Introduction}
We consider a system of $n+1$ hard balls (rods) moving freely on a
straight line, where the only interaction is taking place elastic
collisions between two adjacent balls. In an elastic collision of
two adjacent hard balls, their velocities are redistributed
according to the laws of conservation of energy and momentum. It is
well known that such a dynamical system is isomorphic to a billiard
inside a polyhedral angle. It is also assumed that multiple
collisions, i.e.\ collisions essential between three or more hard
balls, do not occur, corresponding to that if the billiard ball hits
a corner, its further motion is not defined (with some exceptions).

Estimates of the number of collisions in hard ball systems, more
generally, in semi-dispersing billiards has been studied for a long
time because of their importance for proper generalization of the
Boltzmann equation. Sinai \cite{Sin78} proved the existence of a
uniform estimate of the number of collisions of billiard
trajectories in a polyhedral angle. Gal'perin \cite{Gal81} obtained
an explicit estimate for a system of elastic particles (balls of
zero radii) on a line. And Sevryuk \cite{Sev94} gave a uniform
estimate for billiards in a polyhedral angle in terms of a
geometrical characteristic of the angle. The most general results on
this problem were obtained by Burago, Ferleger and Kononenko
\cite{BFK02}. In particular, they provide an explicit estimate
depending only on masses of balls for generalized hard ball systems
on simply connected Riemannian spaces of non-positive sectional
curvature. However, the maximum number of collisions that a hard
ball system may undergo was known only for systems of three
identical balls in Euclidean space of dimension at least 2, besides
for one-dimensional systems of three balls with different masses
which are almost trivial under the billiard approach. Due to Foch,
an example of initial conditions which led three identical hard
balls to four collisions was known, cf.\ Appendix B in \cite{MC00}.
Thurston and Sandri \cite{TS64} discovered a system of three
identical balls suffering four collisions as well. Then Sandri et
al. \cite{SSN64} conjectured that four is the maximum number of
collisions in 1964. A rigorous proof was published until 1993 by
Murphy and Cohen \cite{MC93}.

One of the features of a hard ball system in one dimension is that
the balls always remain the same order on the line. Since we are
only interested in upper bounds of the number of collisions, the
information on length or distance of the system can be completely
ignored for our purpose. It allows us to reduce directly the
original system to an action of a reflection group, generated by $n$
orthogonal reflections, acting on a sphere in $n$-dimensional
Euclidean space $\mathbb{E}^n$, and then to a numbers game. This
natural method reveals the geometric meaning of the constructions of
"relative mass" and "relative velocity" appeared in \cite{MC00}. The
numbers game is intimately related to Coxeter groups, see, for
example, Chapter 4 in \cite{BB05}, although the main interests there
are different from us.

We number the hard balls $0,1,\ldots,n$ in the order of increasing
coordinates and write $m_i$ for the mass of ball $i$. In
\cite{MC00}, Murphy and Cohen showed that for some initial
conditions at least $n(n+1)\over2$ collisions occur and conjectured
that if $m_i\geqslant{m_{i-1}+m_{i+1}\over2},\ i=1,\ldots,n-1$, then
the maximum number of collisions is $n(n+1)\over2$. The purpose of
this paper is to prove the following theorem.
\newtheorem*{main thm}{Main Theorem}
\begin{main thm}
If $m_i\geqslant\sqrt{m_{i-1}m_{i+1}},\ i=1,\ldots,n-1$, then the
maximum number of collisions is $n(n+1)\over2$.
\end{main thm}

\noindent{\it Remark.} Since $\sqrt{ab}\leqslant{a+b\over2}$, for
$a,b>0$, and equality holds if and only if $a=b$, the main theorem
proves the conjecture of Murphy and Cohen with a weaker assumption.
To see why and how the number $n(n+1)\over2$ arises, consider the
simplest case of equal mass: $m_0=m_1=\cdots=m_n$. Suppose $v_{i-1}$
and $v_i$ are the velocities of ball $i-1$ and ball $i$ respectively
$(1\leqslant i\leqslant n)$ before an elastic collision between
them. Then the post-collision velocities $v'_{i-1}$ and $v'_i$ are,
in general, given by
$$
v'_{i-1}=\frac{(m_{i-1}-m_i)v_{i-1}+2m_iv_i}{m_{i-1}+m_i},\quad
v'_i=\frac{2m_{i-1}v_{i-1}+(m_i-m_{i-1})v_i}{m_{i-1}+m_i}.
$$
In the case of equal mass, the two collision balls simply exchange
their velocities. The consequences will become apparent if we
observe changes of the inversion number of the sequence of
velocities $(v_0,v_1,\ldots,v_n)$, which remains constant between
collisions. The inversion number of a sequence of numbers
$\mathbf{q}=(q_0,q_1,\ldots,q_n)$ is defined as the number of its
inversions, that is,
$$\mathrm{inv}\,(\mathbf{q})=\mathrm{card}\,\left\{(i,j)\bigm|i<j,\ q_i>q_j\right\}.$$
It is obvious that
$0\leqslant\mathrm{inv}\,(\mathbf{q})\leqslant{n(n+1)\over2}$ and
$\mathrm{inv}\,(\mathbf{q})=0$ iff $\mathbf{q}$ is an increasing
sequence. When $v_{i-1}>v_i$, a collision between ball $i-1$ and $i$
exchanges the values of the two velocities in the sequence of
velocities so that its inversion number decreases 1. The collisions
then sort the sequence by binary exchanges until the sequence is in
increasing order, after which there can be no more collision.
Therefore, the total number of collisions equals to the inversion
number of the sequence of the initial velocities. In the proof of
the main theorem, we will construct a sequence (depends on the time)
with the similar property: its inversion number remains constant
between collisions and decreases at least 1 in any collision.

\section{Proof of the Main Theorem}
Let $1_i=(\delta_{i0},\delta_{i1},\ldots,\delta_{in})^T\in
\mathbb{E}^{n+1},\ i=0,1,\ldots,n$, where $\delta_{ij}$ is the
Kronecker delta. Write $v_i$ for the velocity of ball $i$ and set
$$
\mathbf{m}=\sum_{j=0}^n\sqrt{m_j}1_j,\quad
\mathbf{v}=\sum_{j=0}^n\sqrt{m_j}v_j1_j\ \in\mathbb{E}^{n+1}.
$$
Then the momentum and energy of the system read
$(\mathbf{m},\mathbf{v})$ and ${1\over2}||\mathbf{v}||^2$
respectively, where $(\cdot,\cdot)$ is the standard scalar product
on $\mathbb{E}^{n+1}$ and $||\cdot||$ is the norm determined by
the scalar product. For $i=1,\ldots,n$, let
$$
\alpha_i=\left(\frac{1_i}{\sqrt{m_i}}-\frac{1_{i-1}}{\sqrt{m_{i-1}}}\right)\bigg/
\left|\left|\frac{1_i}{\sqrt{m_i}}-\frac{1_{i-1}}{\sqrt{m_{i-1}}}
\right|\right|
$$
and $\sigma_i$ be the orthogonal reflection with respect to the
hyperplane passing through the origin with $\alpha_i$ as a unit
normal, that is,
$$\sigma_i:\ \beta\mapsto\beta-2(\alpha_i,\beta)\alpha_i.$$
It is readily seen that
$$
(\alpha_i,\mathbf{m})=0,\ (\alpha_i,\alpha_i)=1,\quad
i=1,\ldots,n,
$$
$$(\alpha_i,\alpha_j)=0,\quad|i-j|>1,$$
$$
(\alpha_i,\alpha_{i+1})=-\frac{1}{m_i}\cdot
\frac{1}{\sqrt{\frac{1}{m_i}+\frac{1}{m_{i-1}}}}\cdot
\frac{1}{\sqrt{\frac{1}{m_{i+1}}+\frac{1}{m_i}}},\quad
i=1,\ldots,n-1,
$$
and $\mathbf{m},\alpha_1,\ldots,\alpha_n$ form a basis of
$\mathbb{E}^{n+1}$ as a vector space. A collision between ball
$i-1$ and ball $i$ is now realized geometrically by the reflection
$\sigma_i$, according to momentum conservation
$(\mathbf{m},\mathbf{v})=(\mathbf{m},\sigma_i(\mathbf{v}))$ and
energy conservation $||\mathbf{v}||^2=||\sigma_i(\mathbf{v})||^2$.
A necessary condition for the collision really taking place is
$v_{i-1}>v_i$, equivalently, $(\alpha_i,\mathbf{v})<0$, since
$$(\alpha_i,\mathbf{v})=-(\alpha_i,\sigma_i({\mathbf{v}}))=
\frac{v_i-v_{i-1}}{\sqrt{\frac{1}{m_i}+\frac{1}{m_{i-1}}}}.$$ If
$|i-j|>1$, then $(\alpha_i,\alpha_j)=0$, i.e., $\sigma_i$ commutes
with $\sigma_j$. It reflects the fact that several binary collisions
may take place simultaneously.

We are now in a position to play a numbers game. Let
$\mathbf{p}=(p_1,\ldots,p_n)$ thought of as a position in the
game. A position $\mathbf{p}$ is called nonnegative if $p_i$ is
nonnegative for all $i=1,\ldots,n$. Choose the weights
$$k_{ij}=-2(\alpha_i,\alpha_j),\quad1\leqslant i,j\leqslant n.$$
Thus
$$k_{ii}=-2,\ k_{ij}=k_{ji},\quad1\leqslant i,j\leqslant n,$$
$$k_{ij}=0,\quad|i-j|>1,$$
and
$$
k_{i,i+1}=\frac{1}{m_i}\sqrt{\frac{2}{\frac{1}{m_i}+\frac{1}{m_{i-1}}}\cdot
\frac{2}{\frac{1}{m_{i+1}}+\frac{1}{m_i}}},\quad i=1,\ldots,n-1.
$$
Moves in the game are defined as follows. A firing of $i$ changes
a position $\mathbf{p}$ by adding $p_ik_{ij}$ to the $j$-th
component of $\mathbf{p}$ for all $j$. More explicitly, a firing
of $i$ changes $\mathbf{p}$ in the following way: switch the sign
of the $i$-th component, add $p_ik_{ij}$ to each adjacent
component $p_j$, and leave all other components unchanged. Such a
move is called negative if $p_i<0$. A negative game is one that is
played with negative moves from a given starting position. The
negative game terminates when it arrives a nonnegative position.

A history of the original hard ball system generates an orbit of
the action of the reflection group generated by
$\sigma_1,\ldots,\sigma_n$, which records the elastic collision
sequence. And the orbit corresponds to a negative play sequence of
the numbers game with the weights $k_{ij}$ by setting
$$\mathbf{p}=(p_1,\ldots,p_n)=((\alpha_1,\mathbf{v}),\ldots,(\alpha_n,\mathbf{v})).$$
We will show that the negative game defined as above must always
terminate in $n(n+1)\over2$ steps no matter what the starting
position is and how it is played.

Let $\mathbf{p}=(p_1,\ldots,p_n)$ be a position in the numbers
game. To avoid analysis case by case, from now on let
$k_{i0}=k_{i,n+1}=p_0=p_{n+1}=p_{n+2}=\cdots=0$ and the same
symbol $\mathbf{p}$ denote the augmented position
$(0,p_1,\ldots,p_n,0,0,\ldots)$. (The values of
$p_0,p_{n+1},p_{n+2},\ldots$ do not change in the whole game.)
Define $q_i=\sum_{j=0}^ip_j,\ i=0,1,2,\ldots,$ and
$\mathbf{q}=(q_0,q_1,\ldots,q_n)$. We will call $\mathbf{q}$ the
potential associated to the position $\mathbf{p}$. Then a position
is nonnegative if and only if its potential is an increasing
sequence. Suppose now we fire $i$, $1\leqslant i\leqslant n$. The
augmented position after the firing is
$$
\mathbf{p}'=(p_0,\ldots,p_{i-2},p_{i-1}+p_ik_{i,i-1},-p_i,p_{i+1}+p_ik_{i,i+1},
p_{i+2},p_{i+3},\ldots),
$$
and hence the potential associated to it becomes
$\mathbf{q}'=(q'_0,q'_1,\ldots,q'_n)$ where
$$
q'_j=
\begin{cases}
q_j,&j\leqslant i-2,\\
q_i-p_i(1-k_{i,i-1}),&j=i-1,\\
q_{i-1}-p_i(1-k_{i,i-1}),&j=i,\\
q_j-p_i(1-k_{i,i-1}+1-k_{i,i+1}),&j\geqslant i+1.
\end{cases}
$$

Using the elementary inequality
$\frac{2}{\frac{1}{a}+\frac{1}{b}}\leqslant\sqrt{ab}$, for
$a,b>0$, we have
$$
k_{i,i+1}\leqslant\frac{1}{m_i}\sqrt{\sqrt{m_im_{i-1}}\cdot\sqrt{m_{i+1}m_i}}
=\sqrt{\frac{\sqrt{m_{i-1}m_{i+1}}}{m_i}},\quad i=1,\ldots,n-1.
$$
If $m_i\geqslant\sqrt{m_{i-1}m_{i+1}}$, then $k_{i,i+1}\leqslant1$,
i.e.\ $(\alpha_i,\alpha_{i+1})\geqslant-{1\over2},\ i=1,\ldots,n-1$.
It follows that, when $p_i<0$, equivalently, $q_{i-1}>q_i$, the
sequence
$$-p_i(0,\ldots,0,1-k_{i,i-1},1-k_{i,i-1},1-k_{i,i-1}+1-k_{i,i+1},\ldots)$$
is increasing. Therefore, the inversion number of the potential
after firing $i\ (1\leqslant i\leqslant n)$
$$
\mathrm{inv}\,(\mathbf{q}')\leqslant\mathrm{inv}\,
(q_0,\ldots,q_{i-2},q_i,q_{i-1},q_{i+1},\ldots,q_n)=\mathrm{inv}\,(\mathbf{q})-1.
$$
The proof is completed since
$0\leqslant\mathrm{inv}\,(\mathbf{q})\leqslant{n(n+1)\over2}$.

\section*{Acknowledgements}
I thank professor Gu,~C.H. for initiating my interests on the
problem and helpful suggestions for developing the manuscript.
Thank also professor Hu,~H.S. for constant support and
encouragement.

\end{document}